\documentclass[11pt]{article}
\usepackage{amssymb}
\usepackage{amsmath}
\usepackage{amsthm}
\usepackage{latexsym}
\usepackage{amsfonts}
\usepackage{graphicx}
\usepackage{graphics}
\usepackage[T1]{pbsi}

\newtheorem{thm}{Theorem}[section]   
\newtheorem{lem}{Lemma}
\newtheorem{cor}{Corollary}
\newtheorem{Def}{Definition}

\theoremstyle{definition}
\newtheorem*{Proof}{Proof}

\newcommand{\dis}{\displaystyle}
\newcommand{\fa}{\forall}

\newcommand{\ra}{\;\rightarrow\;}

\newcommand{\bi}{\beta}

\newcommand{\de}{\delta }

\newcommand{\e}{\varepsilon }
\newcommand{\f}{\varphi}

\newcommand{\vthi}{\vartheta }

\newcommand{\la}{\lambda }
\newcommand{\mi}{\mu }

\newcommand{\si}{\sigma }

\newcommand{\oo}{\omega}

\newcommand{\R}{\mathbb{R}}
\newcommand{\Z}{\mathbb{Z}}

\newcommand{\ssum}{\sum\limits}

\newcommand{\ct}{{\cal{T}}}

\newcommand{\cm}{{\cal{M}}}

\newcommand{\ld}{\ldots}

\newcommand{\sm}{\smallsetminus}

\newcommand{\pr}{\prime\prime}

 \newcommand{\loc}{\mbox{\footnotesize loc}}

\begin{document}

\title{\bf Properties of extremal sequences for the Bellman function of the dyadic maximal operator}
\author{Eleftherios N. Nikolidakis}
\date{}
\maketitle
\noindent
{\bf Abstract:} We prove a necessary condition that has every extremal sequence for the Bellman function of the dyadic maximal operator. This implies the weak-$L^p$ uniqueness for such a sequence.
\section{Introduction} 
\noindent

The dyadic maximal operator on $\R^n$ is defined by
\begin{eqnarray}
\cm_d\phi(x)=\sup\bigg\{\frac{1}{|Q|}\int_Q|\phi(u)u:x\in Q,\;Q\subseteq\R^n \;\text{is a dyadic cube}\bigg\}  \label{eq1.1}
\end{eqnarray}
for every $\phi\in L^1_{\loc}(\R^n)$, where $|\cdot|$ is the Lebesgue measure on $\R^n$ and the dyadic cubes are those formed by the grids $2^{-N}\Z^n$, $N=0,1,2,\ld\,.$

It is well known that it satisfies the following weak type (\ref{eq1.1}) inequality:
\begin{eqnarray}
|\{x\in\R^b:\cm_d\phi(x)\ge\la\}|\le\frac{1}{\la}\int_{\{\cm_d\phi\ge\la\}}
|\phi(u)|du,  \label{eq1.2}
\end{eqnarray}
for every $\phi\in L^1(\R^n)$ and $\la>0$.

From (\ref{eq1.2}) it is easy to prove the following $L^p$-inequality
\begin{eqnarray}
\|\cm_d\phi\|_p\le\frac{p}{p-1}\|\phi\|_p.   \label{eq1.3}
\end{eqnarray}
It is easy to see that (\ref{eq1.2}) is best possible, while (\ref{eq1.3}) is sharp as it can be seen in \cite{5}. (See also \cite{1} and \cite{2} for general martingales).

A way of studying the dyadic maximal operator is to find certain refinements of inequalities satisfied for it.

In this direction the Bellman function of two variables for $p>1$, has been introduced by the following way:
\begin{align}
T_p(f,F)=\sup\bigg\{\frac{1}{|Q|}\int_Q(\cm_d\phi)^p:\phi\ge0,\;&\frac{1}{|Q|}\int_Q
\;\phi(u)du=f, \nonumber\\
&\frac{1}{|Q|}\int_Q\phi^p(u)du=F\bigg\}  \label{eq1.4}
\end{align}
where $Q$ is a fixed dyadic cube on $\R^n$ and $0<f^p\le F$.

The function given in (\ref{eq1.4}) has been explicitely computed. Actually, this is done in a much more general setting of a non-atomic probability measure space $(X,\mi)$ where the dyadic sets are now given in a family of sets $\ct$, called tree, which satisfies conditions similar to those that are satisfied by the dyadic cubes on $[0,1]^n$.

Then the associated dyadic maximal operator $\cm_\ct$ is defined by
\begin{eqnarray}
\cm_\ct\phi(x)=\sup\bigg\{\frac{1}{\mi(I)}\int_I|\phi|d\mi:x\in I\in\ct\bigg\},  \label{eq1.5}
\end{eqnarray}
where $\phi\in L^1(X,\mi)$.

Then the Bellman function (for a given $p>1$) of two variables associated to $\cm_\ct$ is then given by
\begin{eqnarray}
S_p(f,F)=\sup\bigg\{\int_X(\cm_\ct\phi)^pd\mi:\phi\ge0,\;\int_X\phi d\mi=f,\;\int_X\phi^pd\mi=F\bigg\},  \label{eq1.6}
\end{eqnarray}
where $0<f^p\le F$.

In \cite{4}, (\ref{eq1.6}) has been found to be $S_p(f,F)=F\oo_p(f^p/F)^p$ where $\oo_p:[0,1]\ra\big[1,\frac{p}{p-1}\big]$ is the inverse function $H^{-1}_p$ of $H_p$ defined on $\big[1,\frac{p}{p-1}\big]$ by $H_p(z)=-(p-1)z^p+pz^{p-1}$.

As a result the Bellman function is independent of the measure space $(X,\mi)$ and the underlying tree $\ct$.

In this paper we study those sequences of functions: $(\phi_n)_n$, that are extremal for the Bellman function (\ref{eq1.6}). That is $\phi_n:(X,\mi)\ra\R^+$, $n=1,2,\ld$ satisfy $\int_X\phi_nd\mi=f$, $\int_X\phi^p_nd\mi=F$ and
\begin{eqnarray}
\lim_n\int_X(\cm_\ct\phi_n)^pd\mi=F\oo_p(f^p/F)^p.  \label{eq1.7}
\end{eqnarray}
In Section 3 we prove the following
\begin{thm}\label{thm1.1}
Let $\phi_n:(X,\mi)\ra\R^+$ be as above. Then for every $I\in\ct$,
\begin{eqnarray}
\lim_n\frac{1}{\mi(I)}\int_I\phi_nd\mi=f \ \ \text{and} \ \ \lim_n\frac{1}{\mi(I)}\int_I\phi^p_nd\mi=F.  \label{eq1.8}
\end{eqnarray}
Additionally:
\[
\lim_n\frac{1}{\mi(I)}\int_I(\cm_\ct\phi_n)^pd\mi=F\oo_p(f^p/F)^p,
\]
for every $I\in\ct$.  \hfill$\square$
\end{thm}

This gives as an immediate result that extremal functions do not exist for the Bellman function. Another corollary is the weak-$L^p$ uniqueness of such a sequence in all interesting cases. In other words if $(\phi_n)_n$, $(g_n)_n$ are extremal sequences for (\ref{eq1.4}), then $\lim_n\int_Q(\phi_n-g_n)hd\mi=0$, for every $h\in L^p(Q)$, where $\frac{1}{p}+\frac{1}{q}=1$.

Then the following questions arise naturally: \vspace*{0.2cm}\\
{\bf Question 1:} Is an extremal sequence strong $L^p$-unique. By this we mean that we ask the following:

Let $(\phi_n)_n,(g_n)_n$ be extremal sequences. Does it hold that
\[
\dis\lim_n\int_X|\phi_n-g_n|^pd\mi=0\,? \vspace*{-0.3cm}
\]
\hfill$\square$ \medskip

Additionally we mention that (\ref{eq1.8}) is not sufficient for a sequence $(\f_n)_n$ to be extremal. So we may ask the following: \vspace*{0.2cm} \\
{\bf Question 2:} Are there any necessary additional conditions that together with (\ref{eq1.8}) guarantee the extemality of the sequence $(\phi_n)_n$\,?  \hfill$\square$ \medskip

In these questions we hope to answer in the near future.
\section{Extremal sequences} 
\noindent

Let $(X,\mi)$ be a non-atomic probability measure space. We give the following
\begin{Def}\label{def2.1}
A set $\ct$ of measurable subsets of $X$ will be called a tree if the following are satisfied:
\begin{enumerate}
\item[i)] $X\in\ct$ and for every $I\in\ct$, $\mi(I)>0$.
\item[ii)] For every $I\in\ct$ there corresponds a finite or countable subset $C(I)$ of $\ct$ containing at least two elements such that
\begin{itemize}
    \item[(a)] the elements of $C(I)$ are disjoint subsets of $I$
    \item[(b)] $I=\cup C(I)$
\end{itemize}
\item[iii)] $\ct=\dis\bigcup_{m\ge0}\ct_{(m)}$, where $\ct_{(0)}=\{X\}$ and
    \[
    \ct_{(m+1)}=\bigcup_{I\in\ct_{(m)}}C(I).
    \]
\item[iv)] The following holds $\dis\lim_{m\ra\infty}\dis\sup_{I\in\ct_{(m)}}\mi(I)=0$. \hfill$\square$
\end{enumerate}
\end{Def}
\begin{Def}\label{def2.2}
Given a tree $\ct$ we define the maximal operator associated to it as follows:
\[
\cm_\ct\phi(x)=\sup\bigg\{\frac{1}{\mi(I)}\int_I|\phi|d\mi:x\in I\in\ct\bigg\}
\]
for every $\phi\in L^1(X,\mi)$.  \hfill$\square$
\end{Def}

From \cite{4} we obtain the following:
\begin{thm}\label{thm2.1}
The following holds
\[
\sup\bigg\{\int_X(\cm_\ct\phi)^pd\mi:\phi\ge0,\;\int\phi d\mi=f,\;\int_X\phi^pd\mi=F\bigg\}=F\oo_p(f^p/F)^p,
\]
for $0<f^p\le F$. \hfill$\square$
\end{thm}

At last we give the following
\begin{Def}\label{def2.3}
Let $(\phi_n)_n$ be a sequence of non-negative measurable functions defined on $X$ and $0<f^p\le F$, $p>1$. $(\phi_n)_n$ is called $(p,f,F)$ extremal, or simply extremal if the following hold:
\[
\int_X\phi_nd\mi=f,\;\int_X\phi^p_nd\mi=F, \ \ \text{for every} \ \ n=1,2,\ld
\]
\[
\lim_n\int_X(\cm_\ct\phi_n)^pd\mi=F\oo_p(f^p/F)^p.
\]
\end{Def}
\section{Main theorem} 
\begin{thm}\label{thm3.1}
Let $(\phi_n)_n$ be an extremal sequence. Then for every $I\in\ct$ the following hold:
\begin{enumerate}
\item[i)] $\dis\lim_n\frac{1}{\mi(I)}\int_I\phi_nd\mi=f$
\item[ii)] $\dis\lim_n\frac{1}{\mi(I)}\int_I\phi^p_nd\mi=F$
\item[iii)] $\lim\frac{1}{\mi(I)}\int_I(\cm_\ct\phi_n)^pd\mi=F\oo_p(f^p/F)^p$.
\end{enumerate}
\end{thm}
\begin{Proof}
We remind that $\ct_{(0)}=\{X\}$ and $\ct=\dis\bigcup_{m\ge0}\ct_{(m)}$. We prove this theorem for $I\in\ct_{(1)}$. Then inductively it holds for every $I\in\ct_{(m)}$, $m\ge1$.

Suppose then that $\ct_{(1)}=\{I_k,\;k=1,2,\ld\}$ and $I=I_1$. We now set
\[
f^1_n=\frac{1}{\mi(I_1)}\int_{I_1}\phi_nd\mi, \ \ f^2_n=\frac{1}{\mi(X\sm I_1)}\int_{X\sm I_1}\phi_nd\mi,
\]
\setcounter{equation}{0}
\begin{eqnarray}
F^1_n=\frac{1}{\mi(I_1)}\int_{I_1}\phi^p_nd\mi, \ \ F^2_n=\frac{1}{\mi(X\sm I_1)}\int_{X\sm I_1}\phi^p_nd\mi, \ \ \text{for} \ \ n=1,2,\ld,  \label{eq3.1}
\end{eqnarray}

The above sequences are obviously bounded, so passing to a subsequence we may suppose that
\[
\lim_nf^i_n=f^i \ \ \text{and} \ \ \lim_nF^i_n=F^i, \ \ \text{for} \ \ i=1,2.
\]
For any $J\in\ct$ define
\[
\cm_J\phi(t)=\sup\bigg\{\frac{1}{\mi(K)}\int_K|\phi|d\mi:t\in K\in\ct_J\bigg\}, \ \ \text{for} \ \ t\in J,
\]
where $\ct_J$ is defined by
\[
\ct_J=\{K\in\ct:K\subseteq J\}.
\]

Consider the measure space $\big(J,\frac{\mi(\cdot)}{\mi(J)}\big)$, the tree $\ct_J$ and the associated maximal operator $\cm_J$. Then using Theorem \ref{thm1.1}, we have that
\begin{eqnarray}
\frac{1}{\mi(J)}\int_J(\cm_J\phi)^pd\mi\le\frac{1}{\mi(J)}\int_J\phi^pd\mi\cdot
\oo_p\left(\frac{\big(\frac{1}{\mi(J)}\int_J\phi d\mi\big)^p}{\frac{1}{\mi(J)}
\int_J\phi^pd\mi}\right)^p  \label{eq3.2}
\end{eqnarray}
for every $\phi\in L^p(J)$, where $\oo_p:[0,1]\ra\big[1,\frac{p}{p-1I}\big]$ is $H^{-1}_p$, with
\[
H_p(z)=-(p-1)z^p+pz^{p-1}, \ \ z\in \bigg[1,\frac{p}{p-1}\bigg].
\]
(\ref{eq3.2}) now gives (since $H_p$ is decreasing) we have that
\[
H_p\left(\left[\frac{\int_J(\cm_J\phi)^\mi}{\int_J\phi^pd\mi}\right]^{1/p}\right)
\ge\frac{1}{\mi(J)^{p-1}}\frac{\big(\int_J\phi d\mi\big)^p}{\int_J\phi^pd\mi},
\]
which gives
\begin{align}
-&(p-1)\int_J(\cm_J\phi)^pd\mi+p\bigg(\int_J\phi^pd\mi\bigg)^{1/p}\cdot
\bigg(\int_J(\cm_J\phi)^pd\mi\bigg)^{1-\frac{1}{p}}\nonumber \\
&=\frac{1}{\mi(J)^{p-1}}\bigg(\int_J\phi d\mi\bigg)^p+\de_{\phi,J}, \label{eq3.3}
\end{align}
for some $\de_{\phi,J}\ge0$ positive constant depending on $\phi$ and $J$.

For $\phi=\phi_n$ and $J=I_i$, $i=1,2,\ld$ we obtain from (\ref{eq3.3})
\begin{align}
-&(p-1)\int_{I_i}(\cm_{I_i}\phi_n)^pd\mi+p\bigg(\int_{I_i}\phi^p_nd\mi\bigg)^{1/p}\cdot
\bigg(\int_{I_i}(\cm_{I_i}\phi_n)^pd\mi\bigg)^{1-\frac{1}{p}} \nonumber\\
&=\frac{1}{\mi(I_i)^{p-1}}\bigg(\int_{I_i}\phi_nd\mi\bigg)^p+\de_{n,i}, \ \ \text{for every} \ \ n=1,2,\ld \ \ \text{and} \ \ i=1,2,\ld\,.  \label{eq3.4}
\end{align}
Summing relations (\ref{eq3.4}) for $i\ge2$ we obtain
\begin{align}
-&(p-1)\sum^{+\infty}_{i=2}\int_{I_i}(\cm_{I_i}\phi_n)^pd\mi+p\sum^{+\infty}_{i=2}
\bigg(\int_{I_i}\phi^p_nd\mi\bigg)^{1/p}\bigg(\int_{I_i}(\cm_{I_i}\phi_n)^pd\mi
\bigg)^{1-\frac{1}{p}} \nonumber \\
&=\sum^{+\infty}_{i=2}\frac{1}{\mi(I_i)^{p-1}}\bigg(\int_{I_i}\phi_nd\mi\bigg)^p+
\sum^{+\infty}_{i=2}\de_{n,i}.  \label{eq3.5}
\end{align}
In view now of Holder's inequality in its primitive form:
\[
\sum_ia_ib_i\le\bigg(\sum_ia^p_i\bigg)^{1/p}\bigg(\sum_ib^q_i\bigg)^{1/q},
\]
for $a_i,b_i\ge0$ and $q=p/p-1$, (\ref{eq3.5}) gives
\begin{align}
-&(p-1)A_2(n)+p\bigg(\int_{X\sm I_1}\phi^p_nd\mi\bigg)^{1/p}\cdot\big[A_2(n)\big]^{1-\frac{1}{p}}\nonumber \\
&\ge\sum^{+\infty}_{i=2}\frac{1}{\mi(I_i)^{p-1}}\bigg(\int_{I_i}\phi_nd\mi\bigg)^p+
\sum^{+\infty}_{i=2}\de_{n,i}, \ \ \text{where}  \label{eq3.6}
\end{align}
\begin{align}
A_2(n)=\sum^{+\infty}_{i=2}\int_{I_i}(\cm_{I_i}\phi_n)^pd\mi.  \label{eq3.7}
\end{align}
(In the last inequality we used the fact that $X\sm I_1=\dis\bigcup^{+\infty}_{i=2}I_i$).

We use now Holder's inequality in the following form:
\[
\frac{(\la_1+\la_2+\cdots+\la_m)^p}{(\si_1+\si_2+\cdots+\si_m)^{p-1}}
\le\frac{\la^p_1}{\si^{p-1}_1}+
\frac{\la^p_2}{\si^{p-1}_2}+\cdots+\frac{\la^p_m}{\si^{p-1}_m},
\]
where $\si_i$, $\fa\;i=11,2,\ld$ and $\la_i\ge0$, and obtain:
\begin{eqnarray}
\sum^{+\infty}_{i=2}\frac{1}{\mi(I_i)^{p-1}}\bigg(\int_{I_i}\phi_n\mi\bigg)^p\ge
\frac{1}{\mi(X\sm I_1)^{p-1}}\bigg(\int_{X\sm I_1}\phi_nd\mi\bigg)^p=\mi(X\sm I_1)f^2_n.  \label{eq3.8}
\end{eqnarray}
We also set
\begin{eqnarray}
A_3(n)=\int_{X\sm I_1}(\cm_\ct\phi_n)^pd\mi, \ \ \text{for} \ \ n=1,2,\ld\,. \label{eq3.9}
\end{eqnarray}

Then by definition of $\cm_{I_i}$ we have that
\begin{eqnarray}
A_3(n)\ge A_2(n).  \label{eq3.10}
\end{eqnarray}
From the above we have then that:
\begin{eqnarray}
-(p-1)A_2(n)+p\mi(X\sm I_1)^{1/p}(F^2_n)^{1/p}[A_3(n)]^{1-\frac{1}{p}}=\mi(X\sm I_i)(f^2_n)^p+\de^{(1)}_n,  \label{eq3.11}
\end{eqnarray}
where $\de^{(1)}_n\ge\ssum^{+\infty}_{i=2}\de_{n,i}$.

By passing to a subsequence we may suppose that $\dis\lim_nA_3(n)=A_3$.
\end{Proof}

We will use now the following Lemma, the proof of which will be given at the end of this section.
\begin{lem}\label{lem3.1}
If $(\phi_n)_n$ is extremal then we have that
\[
\lim_n\mi(\{\cm_\ct\phi_n=f\})=0. \vspace*{-0.3cm}
\]
\hfill$\square$

\end{lem}
%

From this Lemma and Definitions (3.7) and (3.9) we easily obtain that $\dis\lim_nA_2(n)=\dis\lim_nA_3(n)=A_3$, in view of the fact that $I\in\ct_{(1)}$ for $i=2,3,\ld\,.$ (\ref{eq3.11})  now becomes
\begin{align}
-&(p-1)\int_{X\sm I_1}(\cm_\ct\phi_n)^pd\mi+p\mi(X\sm I_1)^{1/p}(F^2_n)^{1/p}\bigg(\int_{X\sm I_1}(\cm_\ct\phi_n)^pd\mi\bigg)^{1-\frac{1}{p}} \nonumber \\
&=\mi(X\sm I_i)(f^2_n)^p+\de_n^{\prime\prime},    \label{eq3.12}
\end{align}
where $|\de^{\pr}_n-\de^\prime_n|\ra0$, as $n\ra+\infty$.

In the same way we have that:
\begin{align}
-&(p-1)\int_{I_1}(\cm_\ct\phi_n)^pd\mi+p\mi(I_1)^{1/p}(F_n^1)^{1/p}\cdot
\bigg(\int_{I_1}(\cm_\ct\phi_n)^pd\mi\bigg)^{1-\frac{1}{p}} \nonumber \\
&=\mi(I_1)(f^1_n)^p+\e^{\pr}_n,  \label{eq3.13}
\end{align}
where $\e^{\pr}_n$ is such that $|\e^{\pr}_n-\e^\prime_n|\ra0$, $n\ra+\infty$ for some sequence $\e^\prime_n$ for which $\e^\prime_n\ge\de_{n,1}$.

Summing now (\ref{eq3.12}) and (\ref{eq3.13}) and using Holder's inequality in both previously mentioned forms we have that:
\begin{align}
-&(p-1)\int_X(\cm_\ct\phi_n)^pd\mi+pF^{1/p}\bigg(\int_X(\cm_\ct\phi_n)^pd\mi\bigg)^{1-\frac{1}{p}} \nonumber \\
&\ge\mi(I_1)(f^1_n)^p+\mi(X\sm I_1)(f^2_n)^p+\de^{\pr}_n+\e^{\pr}_n\ge f^p+\de^{\pr}_n+\e^{\pr}_n,  \label{eq3.14}
\end{align}
which gives
\begin{eqnarray}
-(p-1)\int_X(\cm_\ct\phi_1)^pd\mi+pF^{1/p}\bigg(\int_X(\cm_\ct\phi_n)^pd\mi\bigg)^{1-\frac{1}{p}}
=f^p+\vthi_n,  \label{eq3.15}
\end{eqnarray}
where $\vthi_n\ge\de^{\pr}_n+\e^{\pr}_n$, $n=1,2,\ld\,.$

The hypothesis now for $(\phi_n)$ is that
\[
\lim_n\int_X(\cm_\ct\phi_n)^pd\mi=F\oo_p(f^p/F)^p.
\]
This gives $\vthi_n\ra0$ in (\ref{eq3.15}), so
\[
\lim_n\de^\prime_n+\lim\e^\prime_n\le0\Rightarrow\de^\prime_n\ra0,\;\e^\prime_n\ra0, \;\de^{\pr}_n\ra0,\;\e^{\pr}_n\ra0.
\]
As a consequence we have
\[
\mi(I_1)(f^1)^p+\mi(X\sm I_1)(f^2)^p=f^p
\]
because of equality in (\ref{eq3.14}), as $n\ra+\infty$.

Since now $\mi(I_1)f^1+\mi(X\sm I_1)f^2=f$ and $t\mapsto t^p$ is strictly convex on $(0,+\infty)$ we have that $f^1=f^2=f$.

Since now $\de^{\pr}_n\ra0$, we have because of (\ref{eq3.12}) and $f^2=f$ that
\begin{eqnarray}
\lim_n\frac{1}{\mi(X\sm I_1)}\int_{X\sm I_1}(\cm_\ct\phi_n)^pd\mi=F_2\oo_p(f^p/F_2)^p.  \label{eq3.16}
\end{eqnarray}
Similarly
\begin{eqnarray}
\lim_n\frac{1}{\mi(I_1)}\int_{I_1}(\cm_\ct\phi_n)^pd\mi=F_1\oo_p(f^p/F_1)^p. \label{eq3.17}
\end{eqnarray}
Since $(\phi_n)_n$ is extremal the last two equations give
\begin{eqnarray}
\mi(I_1)\cdot F_1\oo_p(f^p/F_1)^p+\mi(X\sm I_1)\cdot F_2\oo_p(f^p/F_2)^p=F\oo(f^[/F).  \label{eq3.18}
\end{eqnarray}
But as we shall prove in Lemma \ref{lem3.2} below the following function $t\mapsto t\oo_p(f^p/t)^p$, $t\in(f^p,+\infty)$ is strictly concave. So since $\mi(I_1)F_1+\mi(X\sm I_1)F_2=F$ we have because of (\ref{eq3.18}) that $F_1=F_2=F$ and because of (\ref{eq3.17}):
\[
\lim_n\frac{1}{\mi(I)}\int_I(\cm_\ct\phi_n)^pd\mi=F\oo_p(f^p/F)^p,
\]
and Theorem \ref{thm3.1} is now proved.  \hfill$\square$. \vspace*{0.2cm}
%

We prove now the following
\begin{lem}\label{lem3.2}
Let $G:(1,+\infty)\ra\R^+$ defined by $G(t)=t\oo_p(1/t)^p$. Then $G$ is strictly concave.
\end{lem}
\begin{Proof}
It is known from \cite{4} that $\oo_p$ satisfies
\[
\frac{d}{dx}[\oo_p(x)]^p=-\frac{1}{p-1}\frac{\oo_p(x)}{\oo_p(x)-1}, \ \ x\in[0,1].
\]
So we can easily see that
\[
G^\prime(t)=\oo_p(1/t)^p+\frac{1}{p-1}\frac{1}{t}\frac{\oo_p(1/t)}{\oo_p(1/t)-1}, \ \ \text{and}
\]
\[
G^{\pr}(t)=\frac{1}{p-1}\cdot\frac{1}{t}\bigg(\frac{g(t)}{g(t)-1}\bigg)^\prime,
\]
where $g$ is defined on $(1,+\infty)$ by $g(t)=\oo_p(1/t)$. Since $g^\prime(t)>0$, $\fa\;t>1$, we have that $G^{\pr}(t)<0$, $\fa\;t>1$ and Lemma \ref{lem3.2} is proved. \hfill$\square$
\end{Proof}

We continue now with\vspace*{0.2cm} \\
{\bf Proof of Lemma 3.1:} Let us suppose first that $\phi_n$ are $\ct$-simple functions that is for every $n$, there exists a $m_n$ such that $\phi_n$ is constant on each $I\in\ct_{(m_n)}$. That is $\phi_n$ is $\ct$-good in the sense of \cite{4}, for every $n$. If we look at the proof of Lemma 9 in \cite{4} p. 324-326 we see that in all inequalities (4.20), (4.22), (4.23), (4.24) we should have equality in the limit. So as a result we must have that $\frac{1}{(\bi+1-\bi\rho^n_x)^{p-1}}-\frac{(p-1)\bi\rho^n_x}{(\bi+1)^p}\ra
\frac{1}{(\bi+1)^{p-1}}$, for $\bi=\oo_p(f^p/F)-1$, where $\rho^n_x=\frac{a^n_x}{\mi(x)}=a^n_x$, where $a^n_x=\mi(\{\cm_\ct\phi_n=f\})$. But this can happen only if $a^n_x\ra0$. So the proof is completed in the case of $\ct$-simple functions. As for the general case, it is not difficult to see that if $(\phi_n)_n$ is an extremal sequence of measurable functions, then we can construct a sequence of $\ct$-simple functions such that $\int_Xg_nd\mi=f$, $\int_Xg^p_nd\mi\le F$ and
\[
\lim_n\int_Xg^p_nd\mi=F, \ \ \lim_n\int_X(\cm_\ct\phi_n)^pd\mi=F\oo_p(f^p/F)^p.
\]
Additionally, we can arrange every thing in such a way that $\{\cm_\ct\phi_n=f\}\subseteq\{\cm_\ct g_n=f\}$.

Using the same arguments as before for $(g_n)_n$ we can prove that
$\dis\lim_n\mi(\{\cm_\ct g_n=f\})=0$. So $\dis\lim_n\mi(\{\cm_\ct\phi_n=f\})=0$ and Lemma \ref{lem3.1} is proved.  \hfill$\square$ \vspace*{0.2cm}

We now give some applications of the above.

First we prove the following
\begin{cor}\label{cor3.1}
If $0<f^p<F$ then there do not exist extremal functions for the Bellman function $\ct_p(f,F)$ described in (\ref{eq1.4}).
\end{cor}
\begin{Proof}
Let $\phi$ be an extremal function for (\ref{eq1.4}). Applying Theorem \ref{thm3.1} we see that
\[
\frac{1}{\mi(I)}\int_I\phi d\mi=f \ \ \text{and} \ \ \frac{1}{\mi(I)}\int_I\phi^pd\mi=F,
\]
for every $I$ dyadic subcube of $Q$.

As we can see in \cite{3} inequality (\ref{eq1.2}) implies that the base of dyadic sets of the tree $\ct$ differentiates $L^1(Q)$. That is
\[
\begin{array}{l}
  \phi(x)=f \ \ \text{a.e \ \ and} \\ [0.5ex]
  \phi^p(x)=F \ \ \text{-a.e.}
\end{array}
\]
This gives $f^p=F$, which is a contradiction.  \hfill$\square$
\end{Proof}

We also prove
\begin{cor}\label{cor3.2}
Let $T_p(f,F)$ be described by (\ref{eq1.4}). Then if $(\phi_n)_n$, $(g_n)_n$ are extremal sequences for this function, we must have $\phi_n-g_n\overset{w(L^p)}{\longrightarrow}0$, on $\R^n$ as $n\ra+\infty$.
\end{cor}
\begin{Proof}
Of course we have that
\[
\lim_n\frac{1}{|I|}\int_I\phi_n(u)du=\lim_n\frac{1}{|I|}\int_Ig_n(u)du=f.
\]
So $\dis\lim_n\int_Q(\phi_n-g_n)\xi_I(u)du=0$, for every dyadic subcube $I\subseteq Q$.

Since linear combinations of the characteristic functions of the dyadic subcubes of $Q$ are dense in $L^q(Q)$ we should have that $\dis\lim_n\int_Q(\phi_n-g_n)h=0$, for every $h\in L^q(Q)$, that is $\phi_n-g_n\overset{w(L^p)}{\longrightarrow}0$, as $n\ra+\infty$. \hfill$\square$
\end{Proof}

\end{document}